\documentclass{elsart}
\usepackage{epsfig}
\usepackage{amssymb}

 \usepackage{lineno}

\begin{document}
\begin{frontmatter}


\title{Periodically-forced finite networks of heterogeneous coupled oscillators:
a low-dimensional approach}
\author[L]{Carlo R. Laing\corauthref{cor1}},
\corauth[cor1]{Corresponding author. Ph: +64-9-414-0800 x 41038. Fax: +64-9-441-8136.}
\ead{c.r.laing@massey.ac.nz}
\author[K]{Ioannis G. Kevrekidis}
\ead{yannis@princeton.edu}
%\corauth[cor1]{Corresponding author}
\address[L]{ Institute of Information and Mathematical Sciences,
Massey University,
Private Bag 102-904
North Shore Mail Centre,
Auckland,
New Zealand }
\address[K]{Department of Chemical Engineering and Program in Applied and Computational Mathematics,
Princeton University, Princeton, NJ 08544, USA.}

\begin{abstract}
We study a network of 500 coupled modified van der Pol oscillators. 
The value of a parameter associated with each oscillator
is drawn from a normal distribution, giving a heterogeneous network. 
For strong enough coupling the oscillators all have the same period, and we consider periodic
forcing of the network when it is in this state. 
By exploiting the correlations that quickly develop between the state of an
oscillator and the value of its parameter we obtain an approximate low-dimensional description 
of the system in terms of the first few coefficients in a polynomial chaos expansion. 
Standard bifurcation analysis can then be 
performed on this low-dimensional system, and the results obtained from this predict very well 
the behaviour of the high-dimensional system for any set of realisations of the random parameter.
Situations in which the method begins to fails are also discussed.
\end{abstract}

\end{frontmatter}

%{\bf addresses and in the acknowledgements let's thank SJ and for me 
%let's say the work for IGK was partially supported by DARPA and DOE}

\section{Introduction}
Synchronisation is a common phenomenon in biology and elsewhere~\cite{acebon05,gla01,matmir91}.
It is often studied by investigating the conditions under which
oscillators in a particular network will synchronise~\cite{baeguc91,lai98}.
Periodic forcing of systems is also ubiquitous~\cite{gla01}, 
and so it is natural to study the entrainment of a network
of coupled oscillators by a periodic forcing function.
Many authors have studied small networks of two or three non-identical 
oscillators~\cite{baeguc91}, and larger networks of oscillators that may have
some symmetry~\cite{golste99,lai98} or a particular form of coupling~\cite{renerm00}.
The continuum limit in which there exists an infinite number of oscillators has also been
studied in detail and many results are known for this case~\cite{aristr01,balsas00,strogatz00}. 
However, it is known that finite networks can show behaviour that does not occur in the 
continuum limit~\cite{balsas00,doirin06}. 
In many situations, finite networks are the most realistic way
to model a physical system~\cite{butrin99,doirin06,rubter02}. 
Results for large, finite networks will thus help bridge the gap between
small network dynamics (for which bifurcation analysis is straightforward)
and those for an infinite number of oscillators (where statistical physics
provides the appropriate tools).

In this paper we consider
a large but finite heterogeneous network of coupled oscillators, which are collectively periodically forced.
However, we do not analyse the system exactly; instead we analyse a low-dimensional description
of it.  
This is the ``equation-free'' approach developed by Kevrekdis et al.~\cite{kevgea03}. 
The results here extend those of Moon et al.~\cite{moogha06,mookev06}:
we consider two-variable oscillators, capable of undergoing Hopf bifurcations;
we consider periodic forcing of the network, and we perform bifurcation analysis on the system to
understand how the behaviour of the system changes as parameters are varied. 

\begin{figure}
\leavevmode
\epsfxsize=5.4in
\epsfbox{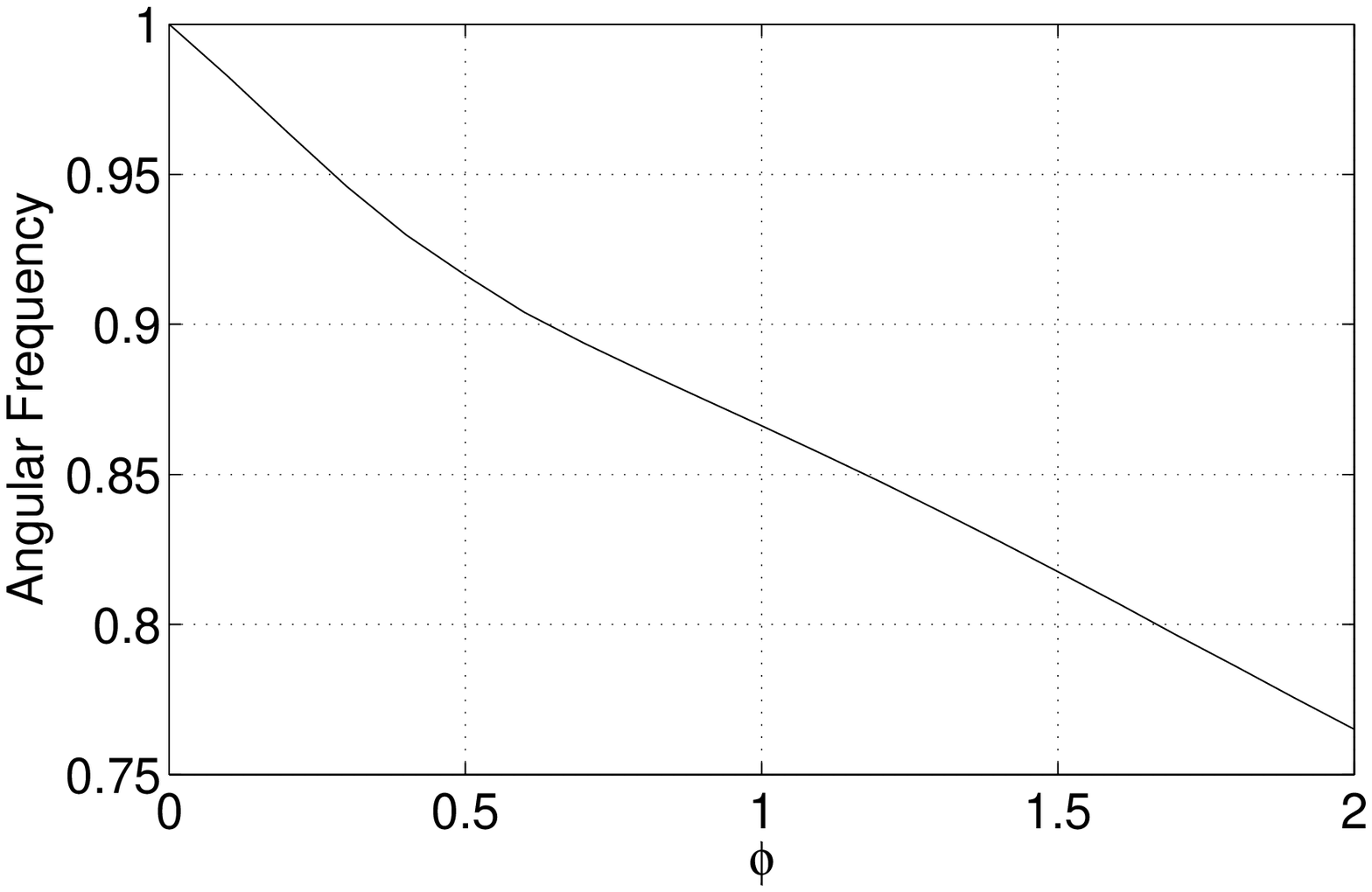}
\caption{Angular frequency of an isolated oscillator ($\beta=0=\epsilon$) as a function of~$\phi$.}
\label{fig:freq1}
\end{figure}

The system we study is
\begin{eqnarray}
   \frac{dx_i}{dt} & = & y_i-x_i\left[x_i^2/3-(\phi+\beta\mu_i)\right]+x_i^2/2-\frac{\epsilon}{N}\sum_{j=1}^{N}(x_i-x_j) \label{eq:dxdt} \\
   \frac{dy_i}{dt} & = & -x_i+A\sin{(\omega t)} \label{eq:dydt}
\end{eqnarray}
for $i=1,\ldots,N$, where $N$ is the number of oscillators in the network. For most of this paper we
set $N=500$. 
The oscillators are van der Pol oscillators~\cite{guchol90} with an extra term $(x_i^2/2)$ which breaks the
internal symmetry $[(x,y)\rightarrow(-x,-y)]$ of each individual oscillator. 
These oscillators were
chosen as being ``typical'' in the sense of not having any particular properties.
For $\beta=0$,
an uncoupled oscillator ($\epsilon=0$)
undergoes a supercritical Hopf bifurcation as $\phi$ increases through zero, with
angular frequency 1. 
The angular frequency for an isolated oscillator
as a function of $\phi$ is shown in Fig.~\ref{fig:freq1}.

The $\mu_i$ are taken from a normal distribution with mean $0$ and standard 
deviation $1$. (As discussed below, the methodology can be used with other distributions.)
If $\beta\neq 0$ the network is heterogeneous, and each oscillator, if uncoupled, would
have a different angular frequency determined by the value of $\phi+\beta\mu_i$. 
When $A=0$, 
for $\beta$ small enough, $\phi$ of moderate size and $\epsilon$ large enough, the oscillators synchronise
in the sense of having the same period. 
Note that oscillators $i$ and $j$ cannot synchronise in the sense of
$x_i(t)=x_j(t)$ for all $t$ unless $\mu_i=\mu_j$. 
In this synchronised state the attractor
of the system is a periodic orbit, which could be parametrised by a periodic variable, say $\theta(t)$.
The variables $x_1,\ldots x_N,y_1,\ldots,y_N$ could each then be written as functions of $\theta$.
This description would no longer be valid if one or more of the oscillators ``unlocked'' from the
group.

We want to study the system in this synchronised state, but do not want to keep track of all
the $2N$ variables $x_1,\ldots x_N,y_1,\ldots,y_N$. 
Instead, we describe the state of the system by a small number of variables. 
We cannot easily derive an equation that governs the dynamics of these
variables, but by repeatedly mapping between the two levels of description of the system we can numerically
evaluate the results of integrating these unavailable equations; we can also find their 
collectively periodic
states and their dependence on parameters, without ever obtaining the reduced equations in closed form.

If the system is in this synchronised state and we increase $A$ from zero,
it will become periodically driven and it may be possible for the oscillators to 
lock with the driving frequency~\cite{gla01}. 
The latter part of this paper will consider this phenomenon in
detail, but we first discuss the particular low-dimensional description of the forced 
system~(\ref{eq:dxdt})-(\ref{eq:dydt}) used here, and how
it can be used in projective integration to speed up direct simulation of this system.

\section{A low-dimensional description}
\begin{figure}
\leavevmode
\epsfxsize=5.2in
\epsfbox{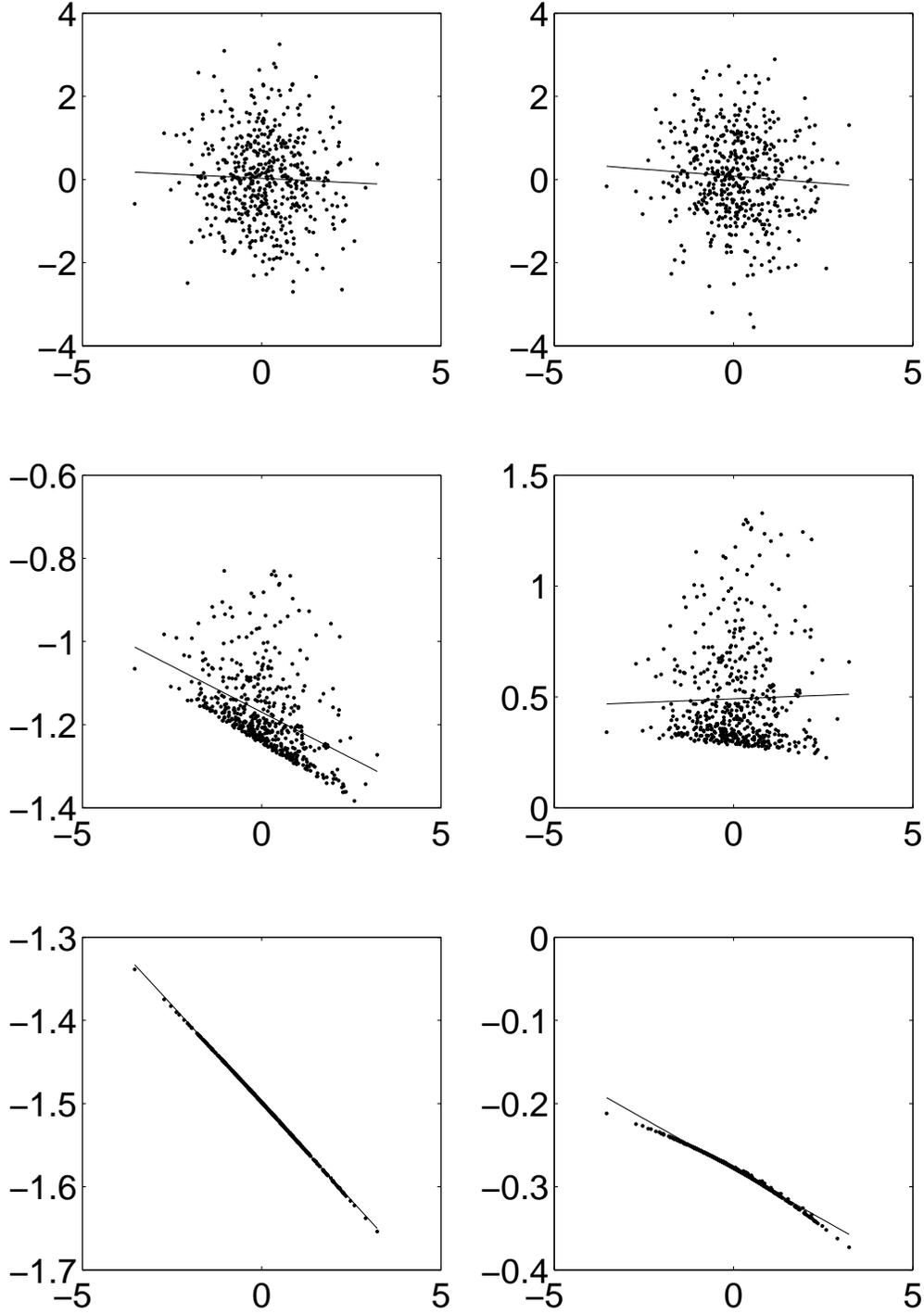}
\caption{Left column: $x_i$ as a function of $\mu_i$ (dots). Right column: $y_i$ as a function 
of $\mu_i$ (dots). Also included are the polynomial chaos expansions with $q=1$ (lines). From
top to bottom: $t=0,2\pi/\omega,4\pi/\omega$. Parameters are 
$A=0.5,\omega=0.85,\beta=0.1,\epsilon=1,N=500$.}
\label{fig:slav}
\end{figure}

The main idea behind the low-dimensional description used here 
depends on correlations that rapidly develop 
between $x_i,y_i$ and the value of $\mu_i$ in the parameter regime where synchronization
eventually prevails.
This is demonstrated in Fig.~\ref{fig:slav} where we plot
the $x_i$ and $y_i$ as functions of $\mu_i$ 
for a particular realisation of the $\mu_i$ at three different times. 
We see that after just two periods of the forcing
strong correlations develop between the state of an oscillator and its $\mu_i$ value. 
We will see that these
correlations occur whether or not the network is synchronised with the forcing.

These correlations allow us to
expand the $x$ and $y$ in certain classes of polynomials of $\mu$~\cite{moogha06};
Hermite polynomials are appropriate for a normal distribution of $\mu$,
while different $\mu$ distributions correspond to different 
polynomial types (the so-called Generalised Polynomial Chaos, or GPC~\cite{xiukar02}). 
We write
\begin{eqnarray}
   x(t,\mu) & = & \sum_{j=0}^q a_j(t)H_j(\mu) \\
   y(t,\mu) & = & \sum_{j=0}^q b_j(t)H_j(\mu)
\end{eqnarray}
where $H_j$ is the $j$th Hermite polynomial [$H_0(x)=1, H_1(x)=2x, H_2(x)=4x^2-2,\ldots$].
This expansion is known as a polynomial chaos expansion~\cite{moogha06}, 
and the $a_i$ and $b_i$ are the
polynomial chaos coefficients.
For a specific realisation of the $\mu_i$, we have $x_i(t)=x(t,\mu_i)$ and similarly for the $y_i(t)$.
Our low-dimensional description then involves the $2(q+1)$ coefficients $a_0,\ldots a_q,b_0,\ldots,b_q$.
This description is approximate, and the approximation becomes better as $q$ is increased.
Given $x_i$ and $y_i$ for a particular set of $\mu_i$, the $a_i$ are
found by minimising the quantity
\begin{equation}
   \sum_{i=1}^N\left[x_i-\sum_{j=0}^q a_jH_j(\mu_i)\right]^2 \label{eq:restricta}
\end{equation}
and the $b_i$ are found by minimising
\begin{equation}
   \sum_{i=1}^N\left[y_i-\sum_{j=0}^q b_jH_j(\mu_i)\right]^2. \label{eq:restrictb}
\end{equation}
This is easily done in Matlab using the ``backslash'' operator to solve an overdetermined linear system.
The operator from the $x_i$ and $y_i$ to the $a_i$ and $b_i$ is referred to as the ``restriction''
operator.
Similarly we construct a ``lifting'' operator:
given the $a_i(t)$ and $b_i(t)$ and a 
particular realisation of the $\mu_i$ we have
\begin{eqnarray}
   x_i(t) & =  & \sum_{j=0}^q a_j(t)H_j(\mu_i) \label{eq:liftx} \\
  y_i(t) & = & \sum_{j=0}^q b_j(t)H_j(\mu_i) \label{eq:lifty}
\end{eqnarray}
Armed with these two operators we can now proceed to numerically solve the
unavailable equation for the polynomial chaos coefficients.

\section{Coarse Projective Integration} \label{sec:proj}
Coarse projective integration entails accelerating the simulation of a high-dimensional system 
by projecting forward in
time using only the variables in a low-dimensional description of the system.
This is accomplished by occasionally
performing short bursts of full simulation of the high-dimensional system in order to 
obtain the numerical information (estimation of the time-derivatives of the
low dimensional description variables) required to perform accurate
projections~\cite{kevgea03,moogha06,mookev06}.
We can use the low-dimensional description in the previous section
for coarse projective integration as follows. 
For convenience,
let the high-dimensional description be the variable
\begin{equation}
    X=[x_1,\ldots,x_N,y_1,\ldots,y_N]\in\mathbb{R}^{2N} \label{eq:X}
\end{equation}
and the low-dimensional polynomial chaos coefficient description be the variable
\begin{equation}
     Z=[a_0,\ldots,a_q,b_0,\ldots,b_q]\in\mathbb{R}^{2(q+1)}. \label{eq:Z}
\end{equation}
Given $X(0)$, 
integrate~(\ref{eq:dxdt})-(\ref{eq:dydt}) forward for $N_1$ steps of size $\delta t$.
Calculate $Z$ at some or all of the times $t=0,\delta t,2\delta t,\ldots,N_1\delta t$ using the
restriction operator. 
Use these values of $Z$ to extrapolate the values of $Z$
to a time $N_2\delta t$ in the future,
i.e.~to time $(N_1+N_2)\delta t$. 
Lift from the value of $Z((N_1+N_2)\delta t)$ to
$X((N_1+N_2)\delta t)$ as detailed above.
Restart the integration of~(\ref{eq:dxdt})-(\ref{eq:dydt}) using
$X((N_1+N_2)\delta t)$ as the initial condition and integrate for a further $N_1$ time steps.
Restrict to $Z$ and repeat the procedure.
If the cost of restricting, extrapolating and lifting is small compared to the cost of integrating
the system~(\ref{eq:dxdt})-(\ref{eq:dydt}) for $N_2$ time steps, this procedure may well
be faster than integrating~(\ref{eq:dxdt})-(\ref{eq:dydt}) directly. 
We expect this to be the case when the full system is characterized by a separation
of time scales; the same principle underpins several analytical reduction techniques
(e.g. centre manifolds, inertial manifolds) but in our case the reduction is obtained
on the fly, from the short full simulation bursts.

\begin{figure}
\leavevmode
\epsfxsize=5.4in
\epsfbox{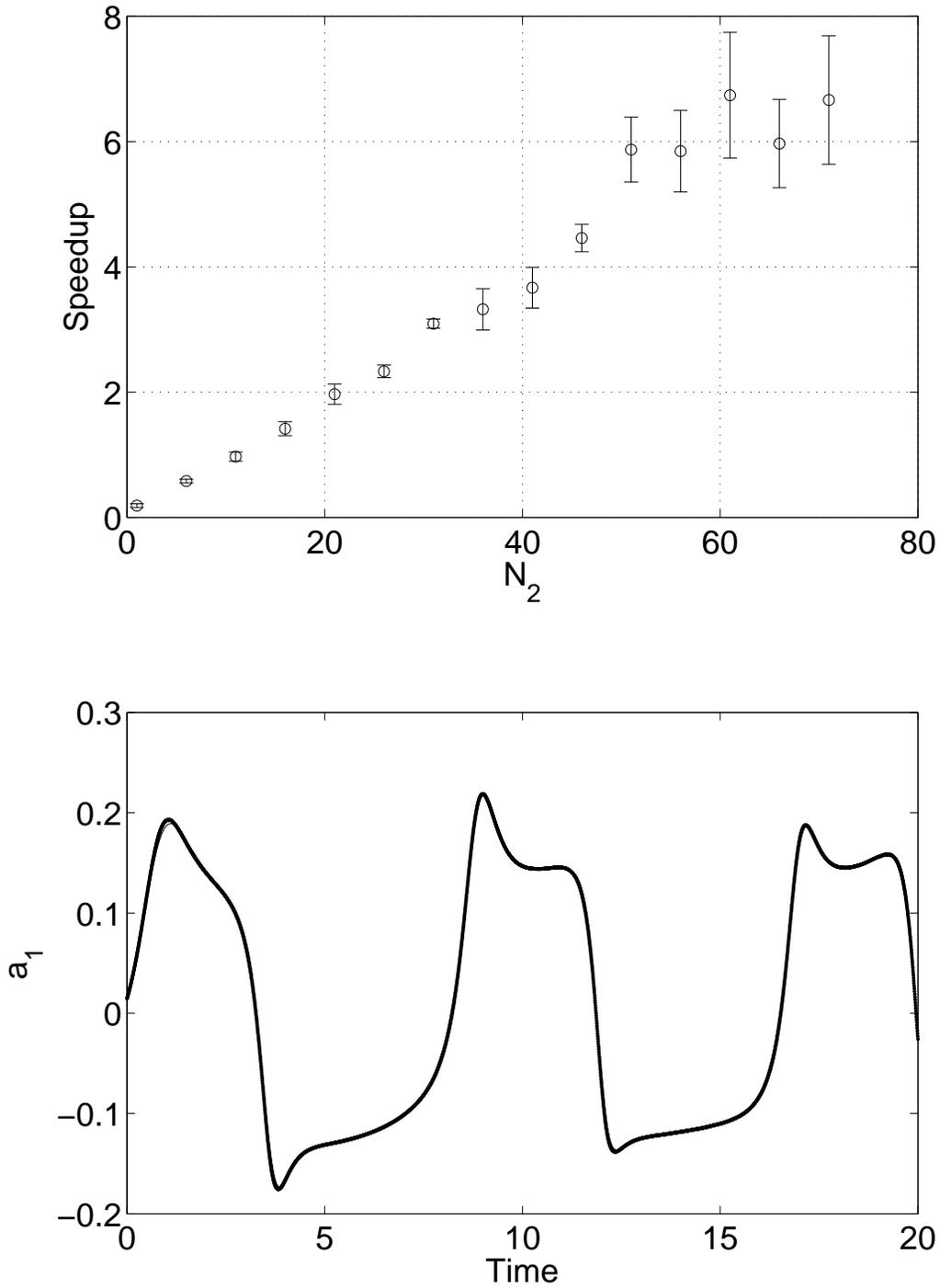}
\caption{Top: speedup (as defined in the text)
as a function of $N_2$, with $N_1=3$. Bottom: $a_1$ for projective
integration with $N_2=1$ (dots) and $a_1$ from full integration (solid line --- indistinguishable from the
dots). 
Other parameters are $A=0.5,\omega=0.85, \phi=1,\beta=0.5,\epsilon=1,q=2$. See text for details.}
\label{fig:proj2}
\end{figure}

\begin{figure}
\leavevmode
\epsfxsize=5.4in
\epsfbox{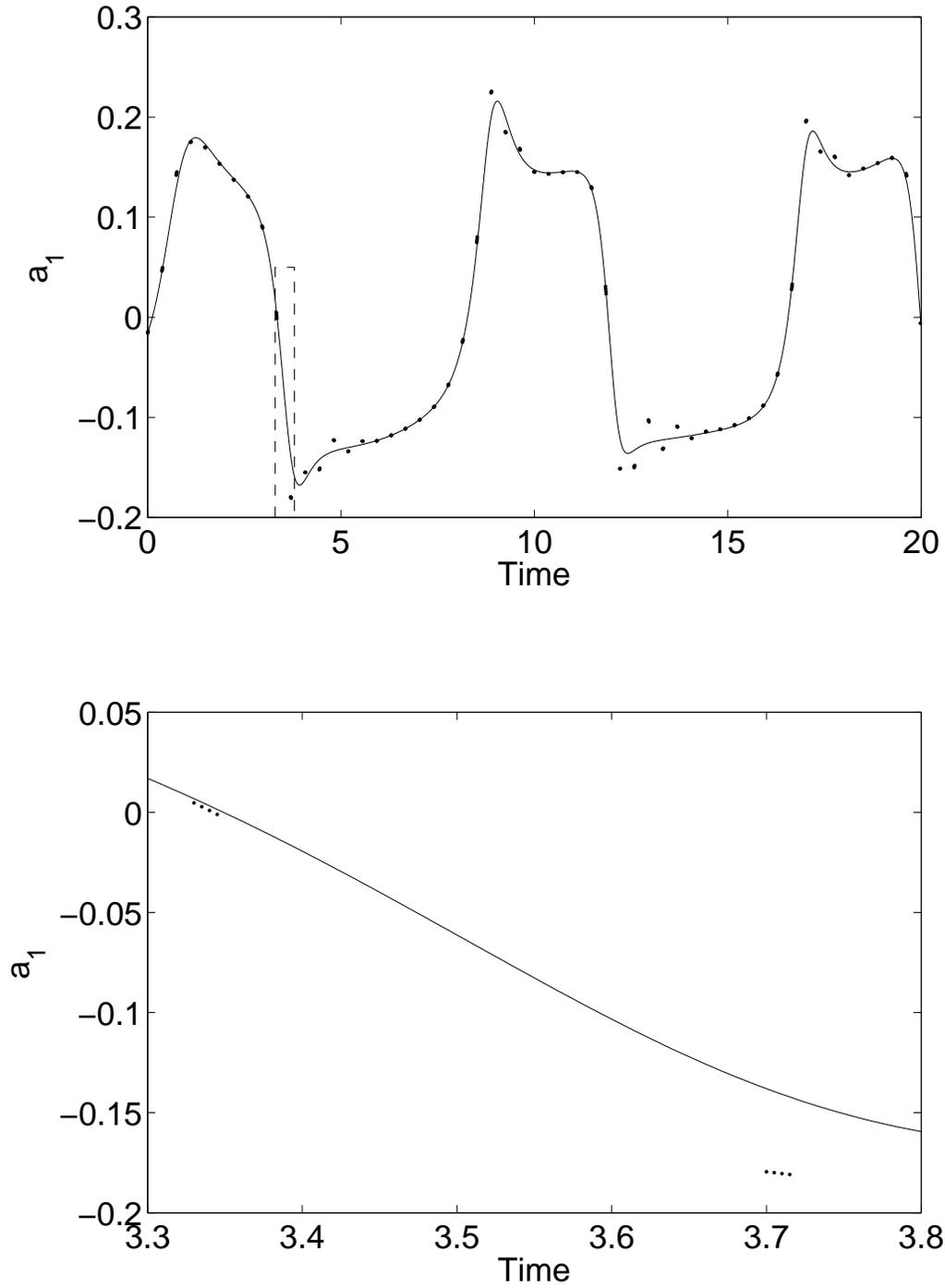}
\caption{Top: $a_1$ for projective
integration with $N_2=71$ (dots) and $a_1$ from full integration (solid line). Bottom:
blowup of the small rectangle in the top panel. Other parameters are 
$A=0.5,\omega=0.85, \phi=1,\beta=0.5,\epsilon=1,q=2$.
See text for details.}
\label{fig:proj1}
\end{figure}

We show results in Figs.~\ref{fig:proj2} and~\ref{fig:proj1} for 
$\delta t=0.005, N_1=3$. 
In this case the projective step involves fitting a cubic to the last $N_1+1$
data points (the last $N_1$ of which were obtained through direct integration of the full system)
and then evaluating this cubic at a time $N_2\delta t$ in the future.
The top panel of Fig.~\ref{fig:proj2} shows the speedup as a 
function of $N_2$. 
The speedup is defined as the time taken to directly
integrate~(\ref{eq:dxdt})-(\ref{eq:dydt})
over $0<t<100$ with time-step $\delta t$ divided by the time taken to integrate over $0<t<100$
using coarse projective integration, as described. 
A speedup greater than 1 ($N_2$ greater than approximately
10) means that projective integration is more efficient than direct integration
(provided that the accuracy in the values of the reduced variables is satisfactory).

Of course, as $N_2$ is
increased the integration will start losing accuracy. 
The bottom panel of Fig.~\ref{fig:proj2}
shows the results of an integration when $N_2=1$ and the two curves shown (one for projective
integration and one for direct integration) are indistinguishable. 
The top
panel of Fig.~\ref{fig:proj1} shows the case when $N_2=71$. 
In this case, coarse projective integration
involves taking 3 steps of length $\delta t$, giving the clusters of 4 points shown in the bottom panel 
of Fig.~\ref{fig:proj1}, and then projecting the $a_i$ and $b_i$ forward a time $71\delta t$,
lifting these values to initialise the $x_i$ and $y_i$ and continuing. 
The integration is clearly less
accurate than that shown in Fig.~\ref{fig:proj2}, but the general behaviour is still qualitatively reproduced.
In the spirit of Taylor series approximations, the extrapolation can only be accurate
up to some fixed interval into the future, so as $\delta t$
is increased, $N_2$ must be decreased, and the speedup will decrease (of course, the accuracy
of the full integration will then also decrease). 
The results shown here will change if $N_1$ is changed
or a different extrapolation scheme is used. 
A full analysis of projective integration for the system
discussed here is beyond the scope of this paper (see discussions 
in~\cite{eeng03,geakev03,kevgea03,samkev06,samroo05};
it is clear, however, that for problems with a large separation 
of time scales and for appropriate parameter choices, it
will be more efficient than straight integration. 
Step adaptation techniques from traditional numerical analysis
based on {\em a posteriori} error estimates can be modified for
the adaptive selection of projective steps.

Note that the simulations shown in Figs.~\ref{fig:proj2} and~\ref{fig:proj1} started at $t=0$,
and thus show transient behaviour, and different initial conditions were used for the two simulations.
It is also important to note that we chose different realisations of the $\mu_i$ at each lifting step;
the results are therefore representative of the expected behaviour {\it over different realizations} of the
random variable. 
Should we only be interested in the acceleration of computations for a particular, single realization,
the results would be even more accurate.

\section{The 1:1 orbit}
Consider the case of 1:1 locking, i.e.~solutions for which each oscillator undergoes one oscillation
during each forcing period. 
The usual way to study this would be to ``strobe'' the system once
each forcing cycle. 
Defining $X_p$ to be the state of the system at $t=2\pi p/\omega$, where $p$ is 
an integer, i.e.
\begin{equation}
   X_p=X(2\pi p)\in\mathbb{R}^{2N}
\end{equation}
where $X$ is defined in~(\ref{eq:X}),
we could construct a map $g:\mathbb{R}^{2N}\rightarrow\mathbb{R}^{2N}$ as 
\begin{equation}
  X_{p+1}=g(X_p)
\end{equation}
A 1:1
locked orbit is then a fixed point of $g$ and its stability is determined by the eigenvalues of the Jacobian of
$g$, evaluated at the fixed point. 
However, finding such a fixed point by, for example, Newton's method,
is computationally very expensive due to the high dimensionality of the system. 
Also, the results we obtain
will only be correct for the particular realisation of the $\mu_i$, a point to which we return below. 
Instead we use
the low-dimensional description of the system in terms of the variable $Z\in\mathbb{R}^{2(q+1)}$.
Defining $Z_p=Z(2\pi p)$, where $Z$ is defined in~(\ref{eq:Z}), we can construct a map
$h:\mathbb{R}^{2(q+1)}\rightarrow\mathbb{R}^{2(q+1)}$ as 
\begin{equation}
   Z_{p+1}=h(Z_p)
\end{equation}
(From now on
we choose $q=1$, so $h:\mathbb{R}^4\rightarrow\mathbb{R}^4$.)

\begin{figure}
\leavevmode
\epsfxsize=5.4in
\epsfbox{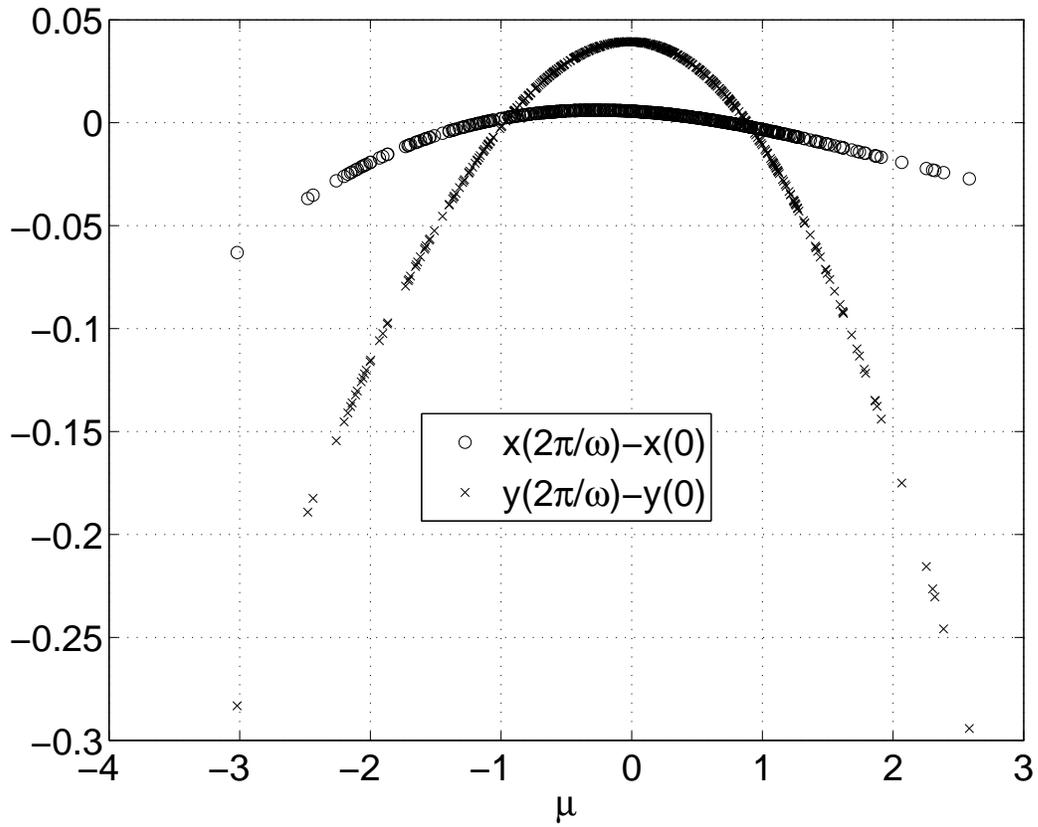}
\caption{Difference in $x$ values after one period (circles), and difference in $y$ values (crosses),
for a fixed point of $h$ (i.e.~for an orbit that is periodic in the coefficients $a_0,\ldots,b_1$).
Parameters are 
$A=0.5,\omega=0.85,\beta=0.5,\epsilon=1$.}
\label{fig:diff}
\end{figure}

Note that a fixed point of $h$ is generally not a fixed point of $g$; 
%%%%%%%%
% when you find a fixed point of the reduced
% the oscillators you found are note a1*H1 ONLY
%
% you started a1*H1 + 0*H2 +0*H3
% and when you finished you were a1*H1 + a2+H2 + a3*H3
% the a2 and a3 YOU do not look at, but it is there, they are FUNCTIONS
% of a1 on an inertial manifold PARAMETRIZED by a1
% and therefore
% a1 is fixed
% but if you want to plot full oscillator description of your fixed point
% you should use the COMPLETE final state
% not just the a1_final * H1
%%%
% we will find the steady state of the inertial form with a1
% not of the TRUNCATION with a1
%
%  maybe worth putting THREE curves in Figure 5, the new of which will
% actually look better -- and we can say what that is 
% 
%%%%%%%%%%%%%%%%%%%%%%%%%%%%%%%%
however, we will see that fixed points
of $h$ do describe the overall behaviour of the system, and the stability follows from the 
eigenvalues of the Jacobian of $h$, evaluated at its fixed points. 
In Fig.~\ref{fig:diff} we show the difference
in $x$ values (and in $y$ values) after a time of one period, for a fixed point of $h$, i.e.~a 1:1
locked orbit in the variables $a_0,\ldots,b_1$. 
We can see that none of the oscillators returned precisely to
its initial condition.
The two distributions $x_i(0)$ and $x_i(2\pi/\omega)$ give the same values
of $a_0$ and $a_1$, even though they clearly do not completely coincide. 
Similarly for $y_i(0)$ and $y_i(2\pi/\omega).$
%
%This is a manifestation of the many-to-one nature of the restriction operator. 
%
If the order of our 
approximation (i.e.~$q$) was increased, the discrepancy shown in Fig.~\ref{fig:diff} would decrease
and the fixed point of $h$ would better approximate the fixed point of $g$.
%%%
%%%Yannis to Carlo
%%% so remember, there is truncation 
%%% and there ALSO is approx inertial form with this
%%% number of variables.

To evaluate $h(Z)$ in practice, we lift from $Z_p$ to $X_p$ using~(\ref{eq:liftx})-(\ref{eq:lifty}), 
integrate~(\ref{eq:dxdt})-(\ref{eq:dydt}) for one period, then restrict from $X_{p+1}$ to $Z_{p+1}$
using~(\ref{eq:restricta})-(\ref{eq:restrictb}). 
Although we could use projective integration
as described in Sec.~\ref{sec:proj}
to integrate~(\ref{eq:dxdt})-(\ref{eq:dydt}), for simplicity here we did not.

\section{Continuation}
\begin{figure}
\leavevmode
\epsfxsize=5.4in
\epsfbox{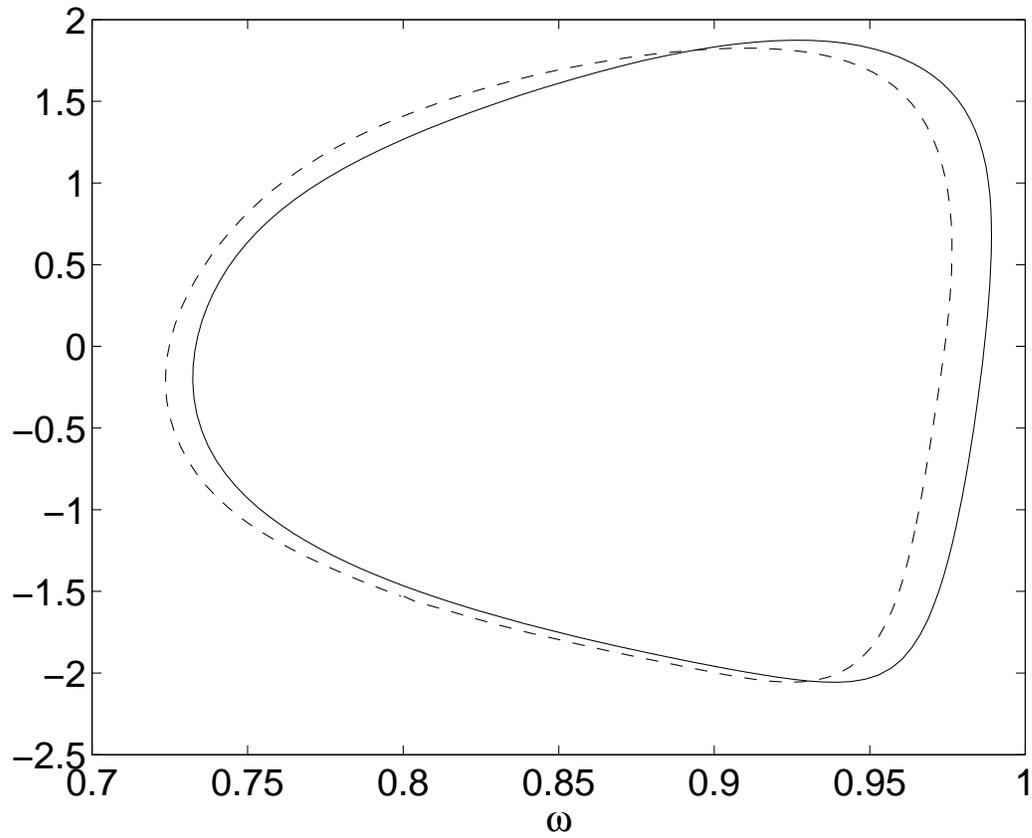}
\caption{The 1:1 orbit. Solid line: $x$ at multiples of $2\pi/\omega$ for a single oscillator.
Dashed line: $a_0$ at multiples of $2\pi/\omega$ for a network with $N=500$ and $\beta=0.5$,
averaged over $r=20$ realisations.
Other parameters are 
$A=0.5,\epsilon=1,\phi=1$.}
\label{fig:lock}
\end{figure}

We can continue fixed points of $h$ as parameters in~(\ref{eq:dxdt})-(\ref{eq:dydt}) 
are varied using standard pseudo-arclength continuation software~\cite{doe97}.
In Fig.~\ref{fig:lock} we show the 
1:1 locked orbit as $\omega$ is varied for a single oscillator (or equivalently, the network
with $\beta=0$, since in this case all oscillators behave identically).
The left and right boundaries of the closed curves are saddle-node bifurcations
where stable and unstable 1:1 locked orbits annihilate one another~\cite{dev89}.

We want to 
analyse the case when $\beta\neq 0$, i.e.~when the network is heterogeneous. 
We could do this for a single
realisation of the $\mu_i$ as above, but to be more general
we choose a number of
different realisations of the $\mu_i$ and average over them. 
We do this averaging within 
our definition of the map $h$. 
Suppose $r$ is the number of realisations we average over.
For each $j=1,\ldots, r$ we calculate $Y_{p+1}^{j}=h(Y_p)$ using the $j$th realisation of the
$\mu_i$. 
(Note that $Y_p$ is fixed.) 
We then define the averaged map $\hat{h}$ as
\begin{equation}
    \hat{h}(Y_p)=\frac{1}{r}\sum_{j=1}^r Y_{p+1}^j. \label{eq:hhat}
\end{equation}
The results of implementing this averaging and following the 1:1 orbit
are shown in Fig.~\ref{fig:lock} (dashed line).
We can see that the effect of the heterogeneity
is to move the range of $\omega$ values for which there is locking to lower 
frequencies. Even though the behaviour of the system
was determined
by following fixed points of $\hat{h}$, the results agree extremely well with those found from direct numerical
integration of the full system~(\ref{eq:dxdt})-(\ref{eq:dydt}) 
for any realisation of the $\mu_i$ from the correct distribution.

\subsection{Varying $A$}
\begin{figure}
\leavevmode
\epsfxsize=5.4in
\epsfbox{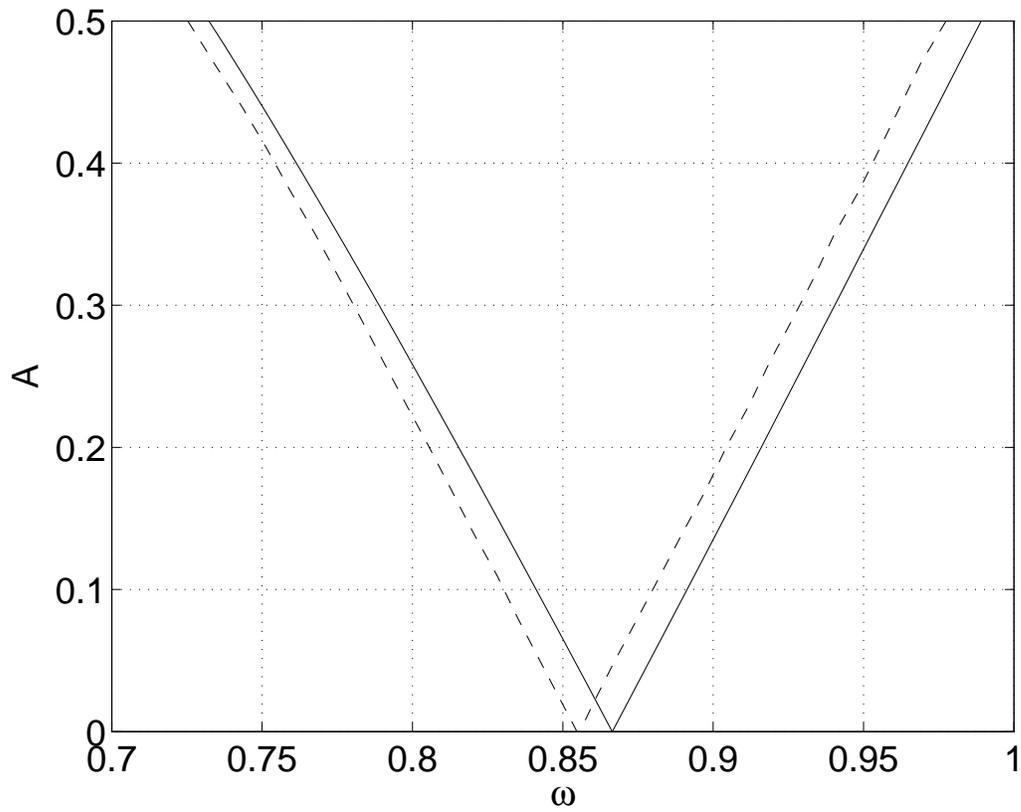}
\caption{Boundaries of the 1:1 Arnold tongue. Solid line: one oscillator, dashed: averaging over $r=20$.
$N=500$ and $\beta=0.5$.
Other parameters are 
$\epsilon=1,\phi=1$.}
\label{fig:one2one}
\end{figure}

We can follow the saddle-node bifurcations of Fig.~\ref{fig:lock}, which mark the edges of the
locking region, as both $A$ and $\omega$ are varied. 
The results are shown in
Fig.~\ref{fig:one2one}, where the resonance ``tongues'' for a single oscillator (solid line) and a network
of 500 oscillators with $\beta=0.5$ (dashed line) are shown. 
Fig.~\ref{fig:lock} is a horizontal slice along the
top boundary of Fig.~\ref{fig:one2one}.

\subsection{Varying $\beta$} \label{sec:varbeta}
\subsubsection{Breakdown of the reduced description}

\begin{figure}
\leavevmode
\epsfxsize=5.4in
\epsfbox{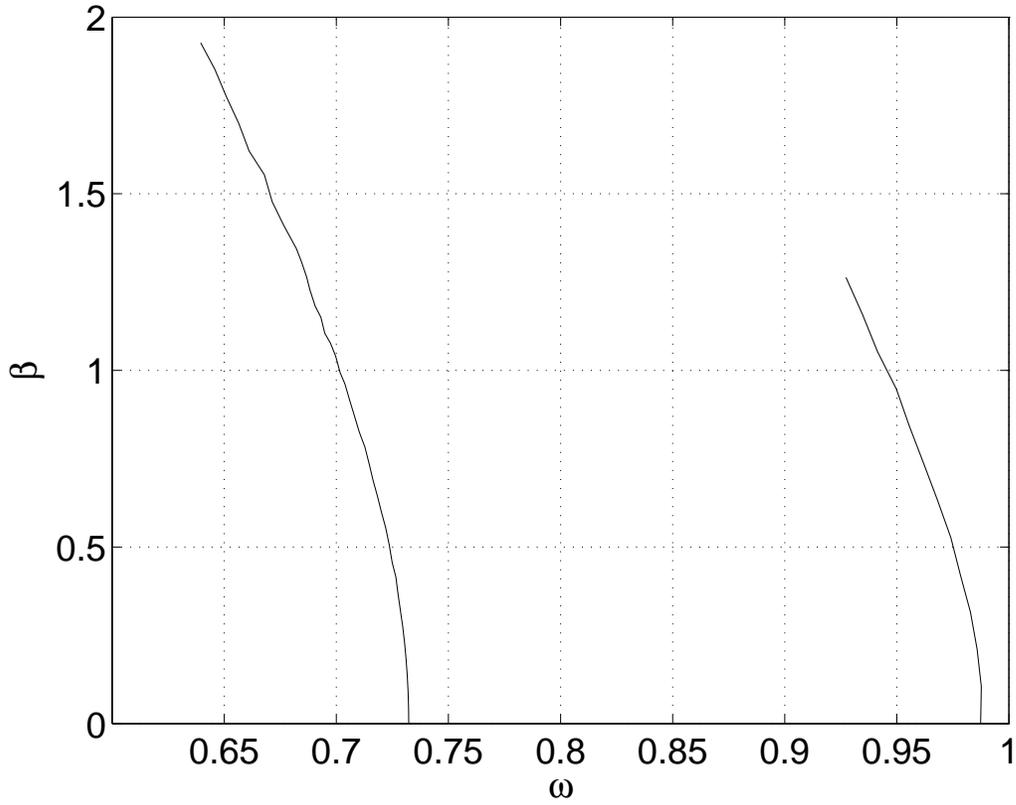}
\caption{Boundaries of the 1:1 Arnold tongue, averaging over between 20 and 60 realisations.
The curves terminate as $\beta$ is increased because the oscillators become too heterogeneous
to synchronise among themselves. $N=500$ and $A=0.5$.
Other parameters are 
$\epsilon=1,\phi=1$.}
\label{fig:betaom}
\end{figure}

\begin{figure}
\leavevmode
\epsfxsize=5.0in
\epsfbox{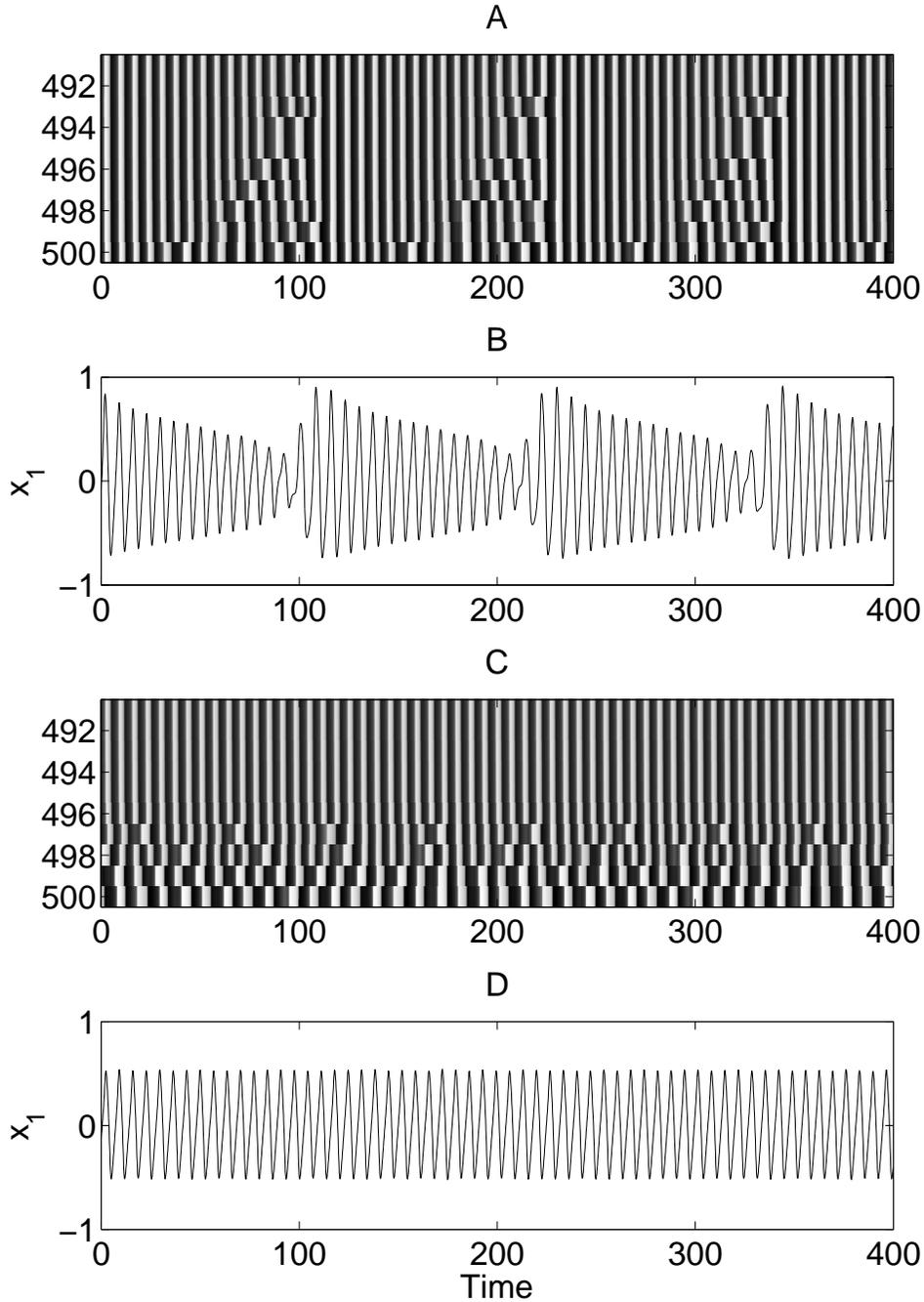}
\caption{A and B: $\omega=0.935,\beta=1.2$ (just to the right of the right boundary shown in 
Fig.~\ref{fig:betaom}). C and D: $\omega=0.925,\beta=1.2$ (just to the left of the right boundary).
Panels A and C show the evolution of 
$x_{491}$ to $x_{500}$ (colour coded, white high), while panels B and D show $x_1$ as a function of time.
$N=500$ and $A=0.5$.
Other parameters are 
$\epsilon=1,\phi=1$.}
\label{fig:edge}
\end{figure}

It is clear that increasing $\beta$ increases the heterogeneity of the network. 
To understand the effects
of this, in Fig.~\ref{fig:betaom} we plot the boundaries of the 1:1 tongue as $\omega$ and $\beta$
are both varied. 
We see the tongue boundaries move to lower frequencies, as expected from 
previous results. 
(Note that Fig.~\ref{fig:lock} shows slices through Fig.~\ref{fig:betaom}
at $\beta=0$ and $\beta=0.5$.) 
When the system is unforced, $\beta$ and $\epsilon$ act 
in opposition: if the
heterogeneity (i.e.~$\beta$) is increased, the coupling strength ($\epsilon$) must be increased in order to keep
the network synchronised. 
However, we consider $\epsilon$ to be fixed. 
Thus for $\beta$ large enough
the forced system will no longer act as a ``super-oscillator'' in which all of the oscillators are synchronised
with each other. 
Once this occurs the concept of locking between all oscillators and the forcing signal
is no longer valid and the algorithm for following
``coarse'' (or macroscopic)  saddle-node bifurcations terminates due to a lack of convergence
within user-specified tolerances.
We demonstrate this phenomenon in detail in Fig.~\ref{fig:edge}.

The top two panels of Fig.~\ref{fig:edge} show the behaviour for a typical realisation of the $\mu_i$
just outside the 1:1 tongue, for a high value of $\beta$. 
The oscillators are ordered by their $\mu$ values.
The behaviour of the 10 oscillators with highest $\mu_i$ is shown in panel A.
In this case, oscillators 1-492 are synchronised with each other, 
but oscillators 493-500 are not synchronised
with the rest of the group. 
However, oscillators 1-492 have also lost their locking to the forcing signal, 
and this is demonstrated
on panel B, where we plot $x_1$ as a function of time. 
This slow (apparently quasiperiodic) modulation is typical for an oscillator
just outside a resonance tongue. 
(Plotting $x_i$ for any $1\leq i\leq 492$ would give a similar picture.)

Panels C and D show the behaviour just inside the tongue (note: a different realisation of the $\mu_i$ from that
in panels A and B has
been used). 
Here, oscillators 1-496 are synchronised with each other, but oscillators 497-500 are
not synchronised with the first 496.
However, now oscillators 1-496 appear  to be still, for all practical purposes, entrained by the forcing.
This is shown in panel D, where $x_1$ is plotted as a function of time. 
The (apparently) periodic oscillation shown here
has the same frequency as the forcing, and a plot of $x_i$ for any $1\leq i\leq 496$ would be very similar.
We say ``apparently'' periodic motion because once one oscillator has desynchronised from the main group,
none of the oscillators will undergo truly periodic motion. Instead, the motion is expected to be
quasiperiodic with at least two frequencies present, or maybe even weakly chaotic.

For the results shown in Fig.~\ref{fig:edge} (with $\beta=1.2$) approximately $1\%$ of the
oscillators (those with the highest values of $\mu$) are not synchronised with the main cluster, either
inside or outside of the tongue. 
However, the remaining $\sim~99\%$ are synchronised with each other
and using the ``macroscopic'' approach taken here we can detect whether this large cluster is synchronised
with the forcing signal or not.

As $\beta$ is increased, the fraction of oscillators no longer locked to the main cluster increases
and the description of the system from the macroscopic point of view as a forced super-oscillator,
using polynomial chaos coefficients,
becomes increasingly flawed. 
This is the reason for deciding to terminate the curves in Fig.~\ref{fig:betaom}.
Note that the two curves in Fig.~\ref{fig:betaom} terminate at different values of $\beta$, but
for both curves, the saddle-node bifurcation following algorithm fails to converge within
tolerances when approximately
$1\%$ of the oscillators become desynchronised from the main group. 

Note that increasing
the number of Hermite polynomials, $q$, used in the macroscopic description
(thus increasing the accuracy of the low-dimensional
description) will not allow
these curves to be followed to greater values of $\beta$. 
It is the lack of synchrony within the
forced network that underlies the termination of the curves. 
Of course, increasing $\epsilon$ would allow the
curves in Fig.~\ref{fig:betaom} to be meaningfully continued to higher values of $\beta$.

Note that if we were to follow a vertical path through the middle of the tongue shown in Fig.~\ref{fig:betaom}
for a particular realisation of the $\mu_i$, there would be many ``fine-scale'' bifurcations as one or more
oscillators desynchronised from the main group. However, these are not visible in our
macroscopic description of the system; we would need to change our macroscopic description in order to detect 
them~\cite{moonab06}.

Moon et al.~\cite{moogha06,mookev06} also considered the loss of synchrony in a heterogeneous network
of Kuramoto oscillators. 
They were studying projective integration and showed that if one or two
oscillators broke from the main cluster, projective integration could continue, as long as the low-dimensional
description was augmented by the phase angle(s) of the oscillators that had lost synchrony. 
%
%Bold et al.~\cite{bolzou06} also considered the loss of synchrony in a network of coupled yeast
%glycolytic oscillators. 
%
%They found that the approach of Moon et al. was also successful in their model,
%but they had to augment their low-dimensional description by six new variables 
%associated with each new ``unlocked" oscillator. 
%
We take a
different approach here, regarding the unsynchronised oscillators as providing a perturbation to the 
dynamics of the synchronised group.

\subsubsection{Phase walkthrough}
For a single periodically driven oscillator, ``phase walkthrough'' can occur just outside a
1:1 resonance tongue~\cite{ermrin84}. In this phenomenon the driven oscillator appears to be
nearly synchronised with the driving oscillator, but every so often 
it undergoes either one extra or one fewer oscillation
than the drive before returning to near synchrony. 
This is because the system lies in the vicinity of a saddle-node
bifurcation of periodic orbits. 
This walkthrough occurs approximately periodically, and
the period scales as $|\omega-\omega^\ast|^{-1/2}$, 
where $\omega^\ast$ is the value of $\omega$ at the
relevant tongue boundary~\cite{doilai02}. 
As can be seen, this slow oscillation can be made arbitrarily slow by adjusting
$\omega$.

A similar phenomenon occurs in our system, but with a slight difference.
For small $\beta$ all of the oscillators are synchronised with one another, effectively acting as one
oscillator, and we can observe phase walkthrough near the tongue boundaries with the scaling 
just mentioned above. 
However, this phenomenon is a result of the system spending a long time in phase space near the remains
of the stable and unstable fixed points of $\hat{h}$, and is thus sensitive to noise or other perturbations. 

Once at least one ``rogue" oscillator has become desynchronised from the rest (as a result of increasing $\beta$)
the system can be thought of as a noisily perturbed oscillator, 
the ``oscillator'' being the vast majority of oscillators
that are synchronised with each other, 
and the ``noise'' resulting from the influence of the desynchronised oscillator(s) on the
rest. 
Thus we expect that we can no longer make the slow oscillation arbitrarily slow just by adjusting $\omega$.
Indeed, for large fixed $\beta$, near the boundaries shown in Fig.~\ref{fig:betaom} 
there is a range of $\omega$ values
for which the slow oscillation (walkthrough) period is not well-defined, 
since perturbations from the desynchronised oscillator(s) affect
the neutrally stable behaviour at the underlying bifurcation, 
resulting in apparently stochastic ``slipping'' relative to the 
forcing signal. It may be possible to describe these rare occurrences in terms of Langevin dynamics on
a low-dimensional
free energy surface~\cite{haasro04,laifre06}.

\subsection{Varying $\phi$}
Another parameter of interest to vary is $\phi$. 
Recall that varying $\phi$ in a single unforced oscillator
causes a Hopf bifurcation, leading to oscillations.
The result of varying $\phi$ is shown in Fig.~\ref{fig:phiom}, for both a single oscillator
and for a network with $\beta=0.3$. 
We see that the tongue terminates at a positive value
of $\phi$, and that heterogeneity moves the tongue boundary
to lower values of $\phi$.

To understand the cusps for low values of $\phi$ we plot in Fig.~\ref{fig:cusp} a cross-section
through Fig.~\ref{fig:phiom} at $\phi=0.8$, for a single oscillator. 
The four saddle-node
bifurcations are clear. 
(For the network, a similar plot is found, not shown.) 
The cusps involve both saddle-node bifurcations being annihilated at a
codimension-two point.  
In the vicinity of these cusps, previous results on the periodically forced van der Pol 
oscillator~\cite{guchol90}
show that there should be a curve of Hopf bifurcations
of the fixed point of the map starting near each cusp, 
which will be associated with the generation of quasiperiodic 
motion.
%
%{\bf look for Peckham in the literature}
%
We can follow these curves using standard algorithms~\cite{doekel91}, and the results 
for the left cusp are shown
in Fig.~\ref{fig:hopf}, both for a single oscillator and for the inhomogeneous network of 500
oscillators. (We also followed the Hopf bifurcation curve associated with the right cusp, not shown.)
The Hopf bifurcations correspond to a complex conjugate pair of eigenvalues crossing
out of the unit circle in the complex plane as $\omega$ is decreased. 
Writing these eigenvalues at bifurcation
as $e^{\pm i\theta}$, we have $\theta=0$ at the rightmost point of the Hopf bifurcation curve (i.e.~eigenvalues
of $+1,+1$) and $\theta$ monotonically increases as $\omega$ is decreased until $\theta=\pi$ (i.e.~eigenvalues
of $-1,-1$) at the leftmost point on the Hopf bifurcation curve.
%At the right extreme of the curve of Hopf bifurcations
%the eigenvalues pass out through 1, i.e.~there is a double eigenvalue of 1. At the leftmost end of the 
%curve of Hopf bifurcations both eigenvalues pass out through $-1$, and as the 
%Hopf bifurcation curve curve is traced from right to left the eigenvalues move around the unit circle
%from $+1,+1$ to $-1,-1$.

Note that while following the curve of Hopf bifurcations for the network, we could not
use a larger value of $\beta$ than $\beta\approx 0.3$ 
(i.e.~we could not use a more heterogeneous network) because
for larger values of $\beta$ the oscillators with the largest values of $\mu$ would become
desynchronised from the rest as the bifurcation was approached. 
The problem
discussed in Sec.~\ref{sec:varbeta} 
regarding the effectiveness of the macroscopic approach would then reoccur.

The Hopf bifurcation for a single forced oscillator is supercritical, with a stable 2--torus being created as the
Hopf bifurcation curve is crossed in the direction of deceasing $\omega$~\cite{guchol90}. 
The criticality of the Hopf
bifurcation for the network seems to be the same as that for a single oscillator, and even though the 
curve in Fig.~\ref{fig:hopf} was found by averaging over 20 realisations of the $\mu_i$, it is still
a very good predictor of the parameter values at which quasiperiodic dynamics occur for any particular 
realisation of the $\mu_i$ (not shown).

Previous results~\cite{guchol90} lead us to
we expect (for each cusp) a curve on which there are orbits homoclinic to the fixed point of the map
(i.e.~homoclinic to a periodic orbit in the full system), emanating from the point where
the curve of Hopf bifurcations and saddle-node bifurcations meet, for both a single oscillator
and the network. We do not consider these curves further. 

%We show the effects of being near this curve in Fig.~\ref{fig:hom}.

%\begin{figure}
%\leavevmode
%\epsfxsize=5.0in
%\epsfbox{hom.eps}
%\caption{Top: $a_0$ for a simulation for one realisation of the $\mu_i.$ Bottom: a Poincare section
%showing $a_0$ and $b_0$ every forcing period (dots) and the fixed point corresponding to the 1:1 orbit
%(circle). As time increases the dots are visited in a clockwise direction.
%Parameters are $\phi=0.71,\omega=0.716,A=0.5,\beta=0.3,\epsilon=1$.}
%\label{fig:hom}
%\end{figure}

%The top panel shows $a_0$ for a simulation very close to the curve of homoclinic bifurcations. The slow
%modulation caused by spending a long time close to a periodic orbit is clear. In the bottom panel
%of Fig.~\ref{fig:hom} we show $a_0$ and $b_0$ at every period of the forcing for the orbit shown
%in the top panel, together with the fixed point (the 1:1 orbit).

\section{Discussion}

In this paper we studied a finite network of heterogeneous, coupled oscillators, all subject to 
the same periodic forcing. We coarse-grained the dynamics, obtaining a low-dimensional description
of the system in terms of a few polynomial chaos coefficients. We defined a return map by sampling
the low-dimensional system once every forcing period; by finding and following fixed points
of this map we
performed standard bifurcation analysis on the 1:1 locked state.
By averaging over realisations of the
distribution of the heterogeneity we have been able to obtain results valid for any particular
realisation.

We have concentrated on only the 1:1 resonance;
the techniques used can be easily applied to any other resonances. 
One issue we have not discussed is varying
$N$, the number of oscillators in the network. 
We found that for small values of $N$
the number of realisations, $r$, of the $\mu_i$ that are averaged over in the definition of 
$\hat{h}$~(eqn.~(\ref{eq:hhat}))
had to be increased in order for continuation algorithms to converge to a given tolerance.
This makes sense since, as $N$ is increased, the difference between a simulation with one particular
realisation of the $\mu_i$ and that for a simulation with a different realisation, will decrease, and thus
fewer realisations will need to be averaged over.

\begin{figure}
\leavevmode
\epsfxsize=5.4in
\epsfbox{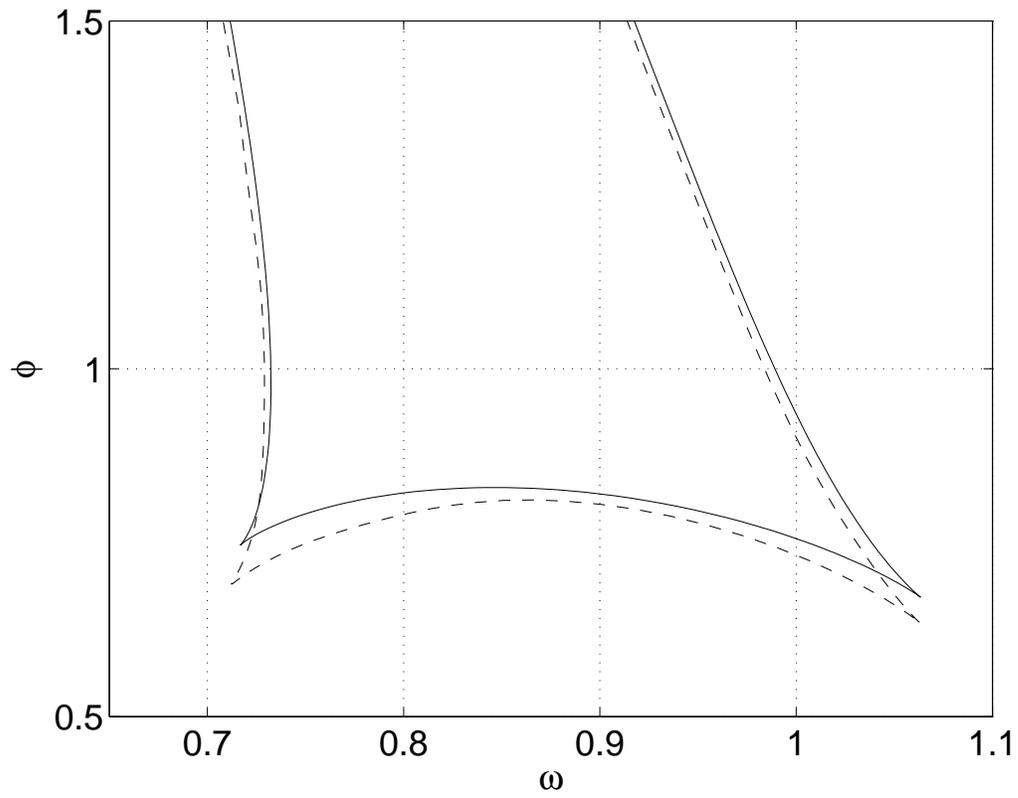}
\caption{Boundaries of the 1:1 orbit. Solid line: one oscillator. Dashed line: a network with
$N=500$ and $A=0.5,\beta=0.3, r=20,\epsilon=1$.}
\label{fig:phiom}
\end{figure}

\begin{figure}
\leavevmode
\epsfxsize=5.4in
\epsfbox{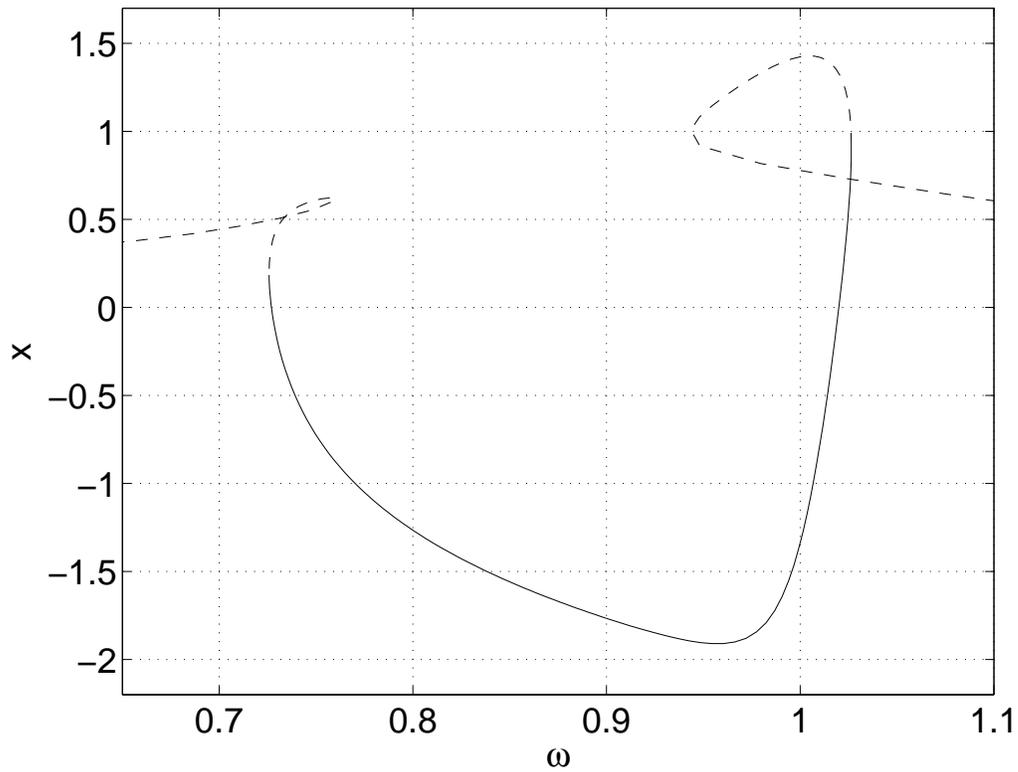}
\caption{$x$ at multiples of $2\pi/\omega$ for the 1:1 orbit of a single oscillator, as $\omega$ is varied.
Solid line: stable, dashed: unstable.
Parameters are $A=0.5,\phi=0.8$. This figure is a horizontal slice at $\phi=0.8$ through
Fig.~\ref{fig:phiom}.}
\label{fig:cusp}
\end{figure}

\begin{figure}
\leavevmode
\epsfxsize=5.4in
\epsfbox{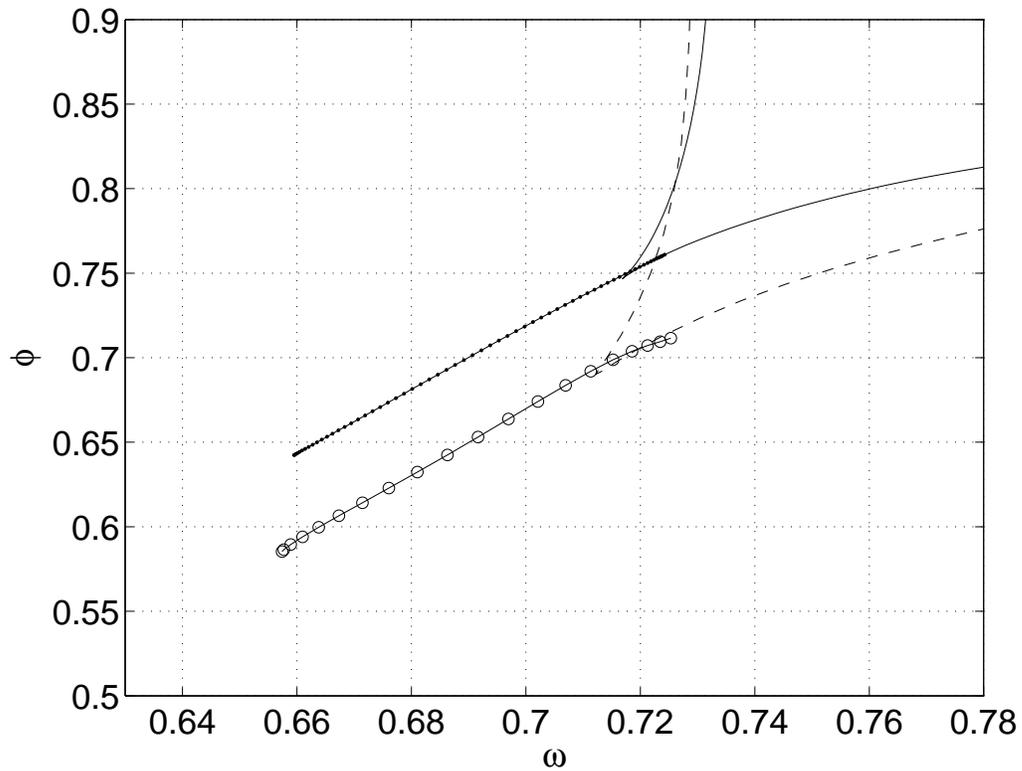}
\caption{Boundaries of the 1:1 Arnold tongue (solid and dashed lines, for a single oscillator and a network
of $N=500$ oscillators, respectively) and Hopf bifurcation curves
(dots and circles joined by lines, for a single oscillator and a network
of 500 oscillators, respectively).
Parameters are $A=0.5,\beta=0.3,r=20,\epsilon=1$. Compare with Fig.~\ref{fig:phiom}.}
\label{fig:hopf}
\end{figure}

Compared to the Kuramoto model results of Moon et al.~\cite{moogha06,mookev06}, or their animal
flocking models \cite{moonab06} where each
oscillator consists of a single variable (a phase angle) the coupled units here are representative
of general ODE-based oscillators, capable
of undergoing Hopf bifurcations. 
It is this extra feature that gives rise to results such as those
in Figs.~\ref{fig:phiom} and~\ref{fig:hopf}. 
%
%Bold et al.~\cite{bolzou06} considered more complicated oscillator
%models (each of which had six variables) but they did not perform any bifurcation analysis. 
%Moon et al.~\cite{moonab06}
%studied what is effectively a network of coupled Kuramoto oscillators in the context 
%of modelling flocking in animals

It would be interesting to apply these ideas to
networks of oscillators where the onset of oscillation is through different types of bifurcations, 
e.g.~homoclinic
bifurcation, saddle-node-on-a-circle (or SNIPER), or a saddle-node bifurcation of periodic orbits, or to other
finite heterogeneous networks that have previously been studied~\cite{butrin99,rubter02}. 
Another
possibility is to study the periodic forcing of networks of coupled bursting neurons. 
The response of isolated
bursting neuron models to periodic forcing has recently been studied~\cite{cooowe01,laicoo05},
as has the behaviour of coupled bursters~\cite{izh00}.

Regarding the problems caused by one or more oscillators desynchronising from the main group
(see Sec.~\ref{sec:varbeta}), one way to deal with this might be to expand the states of oscillators
in functions other than ``globally'' defined polynomials of $\mu$. For example, a wavelet basis
of localised functions may be more suitable~\cite{lemkni04}, 
particularly if the network breaks into clusters, with the oscillators within each
cluster being synchronised.

{\bf Acknowledgements:} We thank Sung Joon Moon for useful conversations about the work presented
here. The work of IGK was partially supported by DARPA and DOE.

\newpage
%\listoffigures

\newpage

\end{document}